\RequirePackage{fix-cm}
\documentclass[twocolumn]{svjour3}          
\smartqed  
\usepackage{graphicx}
\usepackage{amsmath}
\usepackage{amsfonts}
\usepackage{amssymb }
\usepackage{float}
\usepackage{color, colortbl}
\usepackage{array}
\usepackage{tabularx}

\begin{document}

\title{Graphical exploration of the connectivity sets of alternated Julia sets}

\subtitle{$\mathcal{M}$, the set of disconnected alternated Julia sets}

\titlerunning{Graphical exploration of alternated Julia sets}        

\author{Marius-F. Danca\and
Paul Bourke \and
Miguel Romera
}

\institute{Marius-F. Danca \at
Department of Mathematics and Computer Science\\
Avram Iancu University, 400380 Cluj-Napoca, Romania \\
and Romanian Institute for Science and Technology\\
400487 Cluj-Napoca, Romania \\
              \email{danca@rist.ro}
           \and           Paul Bourke\at
The University of Western Australia, Perth, Australia\\
\email{paul.bourke@uwa.edu.au}
\and Miguel Romera \at
Information Security Institute, CSIC,
Serrano 144, 28006 Madrid, Spain\\
\email{miguel@iec.csic.es}
}

\date{Received: date / Accepted: date}

\maketitle

\begin{abstract}
 Using computer graphics and visualization algorithms, we extend in this work the results obtained analytically in \cite{dan1}, on the connectivity domains of alternated Julia sets, defined by switching the dynamics of two quadratic Julia sets. As proved in \cite{dan1}, the alternated Julia sets exhibit, as for polynomials of degree greater than two, the disconnectivity property in addition to the known dichotomy property (connectedness and totally disconnectedness) which characterizes the standard Julia sets. Via experimental mathematics, we unveil these connectivity domains which are four-dimensional fractals. The computer graphics results show here, without substituting the proof but serving as a research guide, that for the alternated Julia sets, the Mandelbrot set consists of the set of all parameter values, for which each alternated Julia set is not only connected, but also disconnected.
\end{abstract}

\section{Introduction}
\label{intro}

In \cite{dan1} the connectedness properties of alternated Julia sets have been analytically and numerically studied. The alternated Julia sets are obtained by alternating two complex quadratic maps

\[{z_{n + 1}} = z_n^2 + {c_i},\,\,\,\,i = 1,2,\,\,\,{z_0},{c_i} \in  \mathbb{C},{\rm{ }}\,~~n \in N{\rm{.}}\]

The alternate iteration of two real Mandelbrot maps was firstly studied in \cite{alme}, and the alternated iteration of two real logistic maps was treated in \cite{dan2} using computer tools such as histograms and surrogate tests. In the case of some classes of continuous and discontinuous dynamical systems of integer and fractional order \cite{dan3}, by switching the parameter value in the underlying initial value problem while it is numerically integrated, any attractor may be synthesized. Moreover, in this way, the chaos control and anti-control of chaos (chaoticization) can be explained and understood in a different way.

As known, the filled Julia set of the complex quadratic polynomial ${P_c}:\mathbb{C} \to \mathbb{C},\,\,\,{P_c}(z) = {z^2} + c,\,\,c \in  \mathbb{C}$, is the set ${K_c} = \left\{ {z \in  \mathbb{C}|P_c^{ \circ n}(z)\not  \to \infty {\rm{ ~~as~~ }}n \to \infty } \right\}$, where $P_c^{ \circ n}$ denotes the $n^{th}$ iteration of $P_c$, while the Julia set of $P_c$ is ${J_c} = \partial {K_c}$, i.e. the boundary of the filled Julia set (see e.g. \cite{peit}).

\begin{remark} \label{rem1}

As it is well known, $J_c$ is connected if and only if $K_c$ is connected. If $c$ belongs to the Mandelbrot set, then the Julia set is connected while, if it is outside the Mandelbrot set, the Julia set is totally disconnected \cite{man}. This was the key result that inspired Mandelbrot, in 1977, to visualize a set in the complex parameter space consisting of all $c$-values that have connected Julia sets. This
set is now known as the Mandelbrot set ($\mathcal{M}$).
\end{remark}

\begin{figure*}
\begin{center}
  \includegraphics[clip,width=0.7\textwidth] {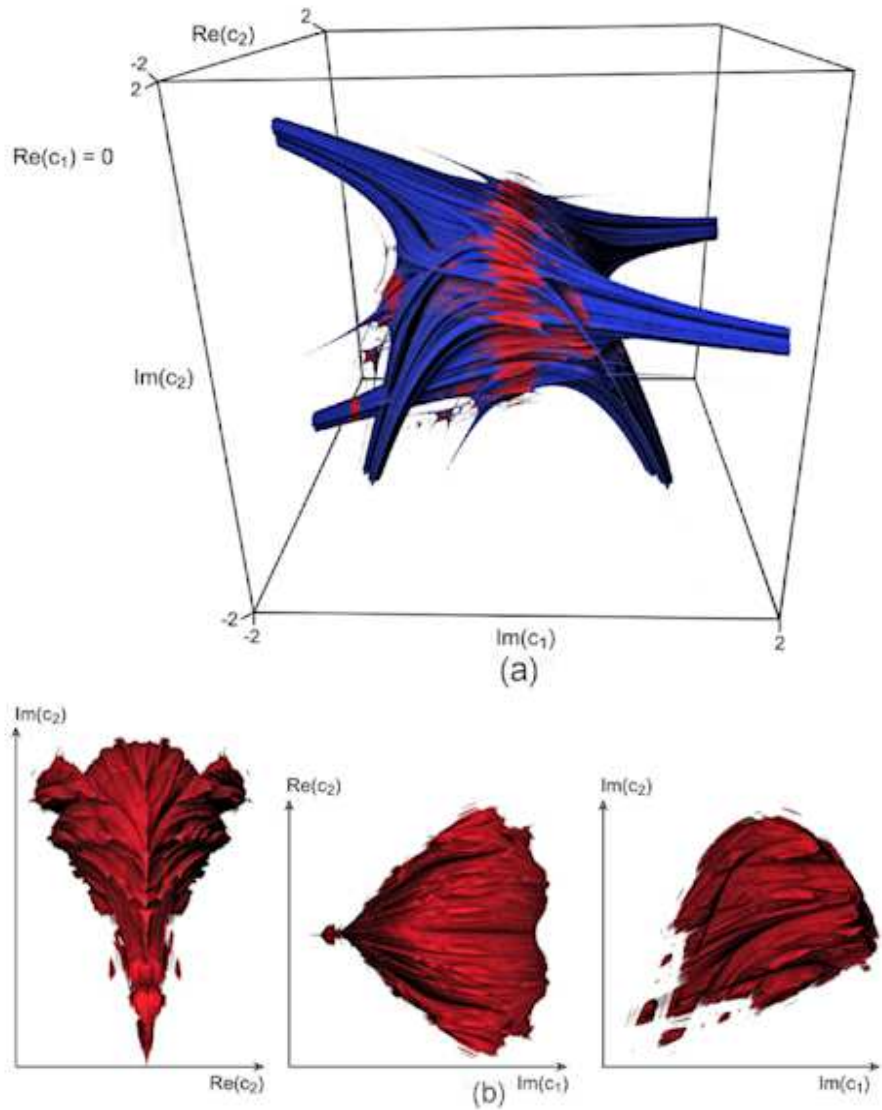}
\caption{The three dimensional body resulting from the hyperplane slice ${\mathop{\rm Re}\nolimits} ({c_1}) = 0$.
 (a) Perspective view of the connectivity (red) and disconnectivity (blue) sets; (b) Projections of the connectivity set along the three orthogonal axes  ${\mathop{\rm Im}\nolimits} ({c_1})$, ${\mathop{\rm Re}\nolimits} ({c_2})$ and ${\mathop{\rm Im}\nolimits} ({c_2})$}
\label{fig1}       
\end{center}
\end{figure*}

\begin{figure*}
\begin{center}
  \includegraphics[clip,width=0.65\textwidth]{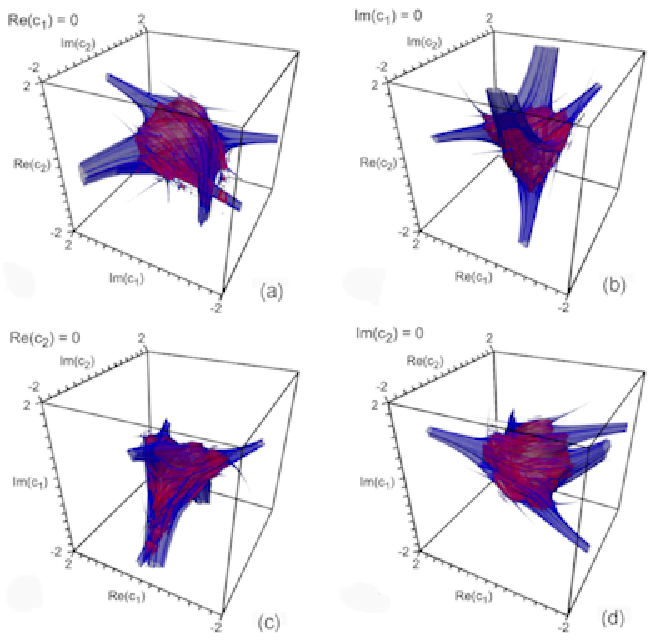}
\caption{Semi-transparent perspective views of $B$ corresponding respectively to: (a) ${\mathop{\rm Re}\nolimits} ({c_1}) = 0$. (b) ${\mathop{\rm Im}\nolimits} ({c_1}) = 0$. (c) ${\mathop{\rm Re}\nolimits} ({c_2}) = 0$. (d) ${\mathop{\rm Im}\nolimits} ({c_2}) = 0$}
\label{fig2}       
\end{center}
\end{figure*}

In order to see if a filled Julia set is connected or totally disconnected\footnote{Intuitively, a \emph{totally disconnected} set is a set which can be be broken into two pieces at each of its points, the breakpoint being always ``in between'' the original set while if a set can be separated into two open and disjoint sets such that neither set is empty and both sets combined give the original set, then the set is called \emph{disconnected}. A set which is not disconnected is \emph{connected}.}  the fate of the orbit of the critical points (where derivative vanishes) must be analyzed. This leads to the well-known fundamental topological dichotomy \cite{fat,jul} for the Julia sets of quadratic polynomials:

\begin{itemize}
\item The Julia set is connected if and only if all the critical orbits are bounded.
\item The Julia set is totally disconnected if (but not only if) all the critical orbits are unbounded. One of the variants of this theorem states that: If the orbit of the critical point is unbounded, then the Julia set is a totally disconnected set (the typical example is the Cantor set, where if the orbit of the critical point is bounded then the Julia set is a connected set). For higher degree polynomials the fundamental dichotomy theorem is no longer valid. Blanchard, Devaney and Keen \cite{dev} introduced the statement for higher degree polynomials. The Julia set is disconnected if at least one critical point iterates to infinity, but it may or may not be totally disconnected.
\item  For higher degree polynomials the fundamental
theorem is no longer valid. Blanchard, Devaney, and Keen introduced
the statement for higher degree polynomials. In this case, the Julia set is
disconnected if at least one critical point iterates to infinity, but it may
or may not be totally disconnected.
\end{itemize}

\begin{figure*}
\begin{center}
  \includegraphics[clip,width=0.75\textwidth]{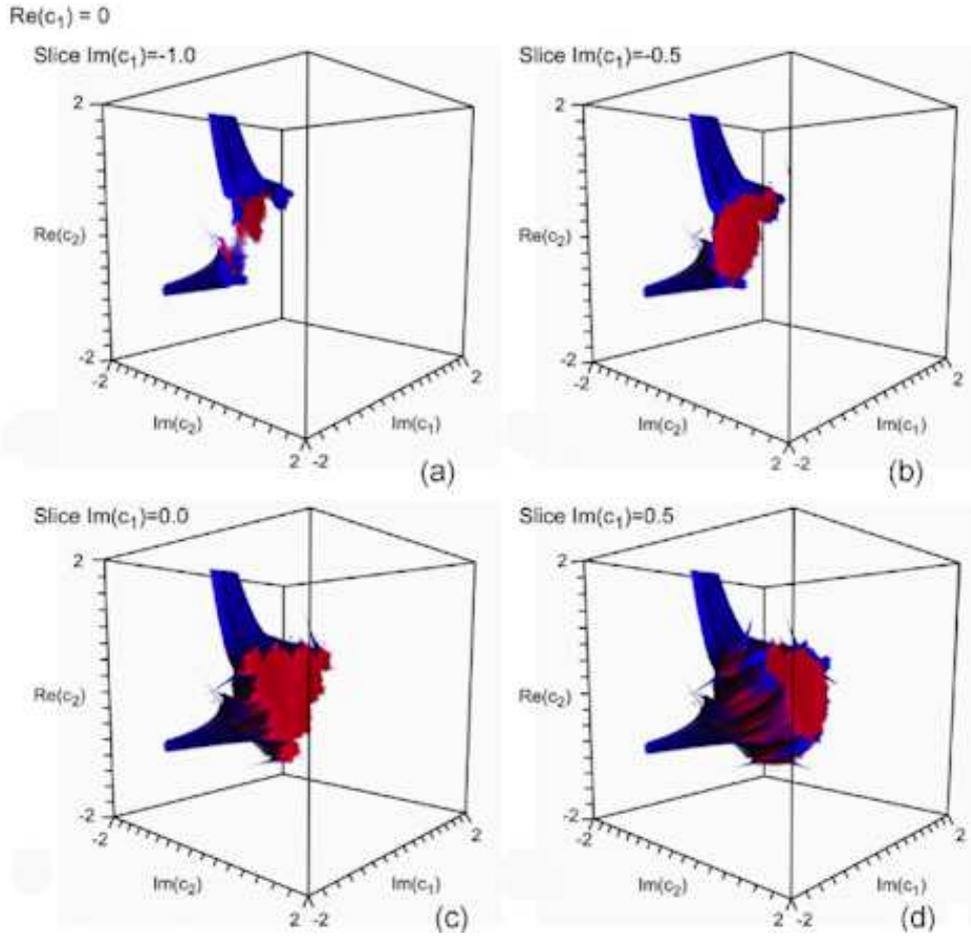}
\caption{Cross sectional views by taking planes of constant $Im(c_2)$ through $B$ for the case $Re(c_1)=0$. (a) $Im(c_1)=-1.0$. (b) $Im(c_1)=-0.5$. (c) $Im(c_1)=0.0$. (d) $Im(c_1)=0.5$}
\label{fig3}       
\end{center}
\end{figure*}

\begin{figure*}
\begin{center}
  \includegraphics[clip,width=0.65\textwidth]{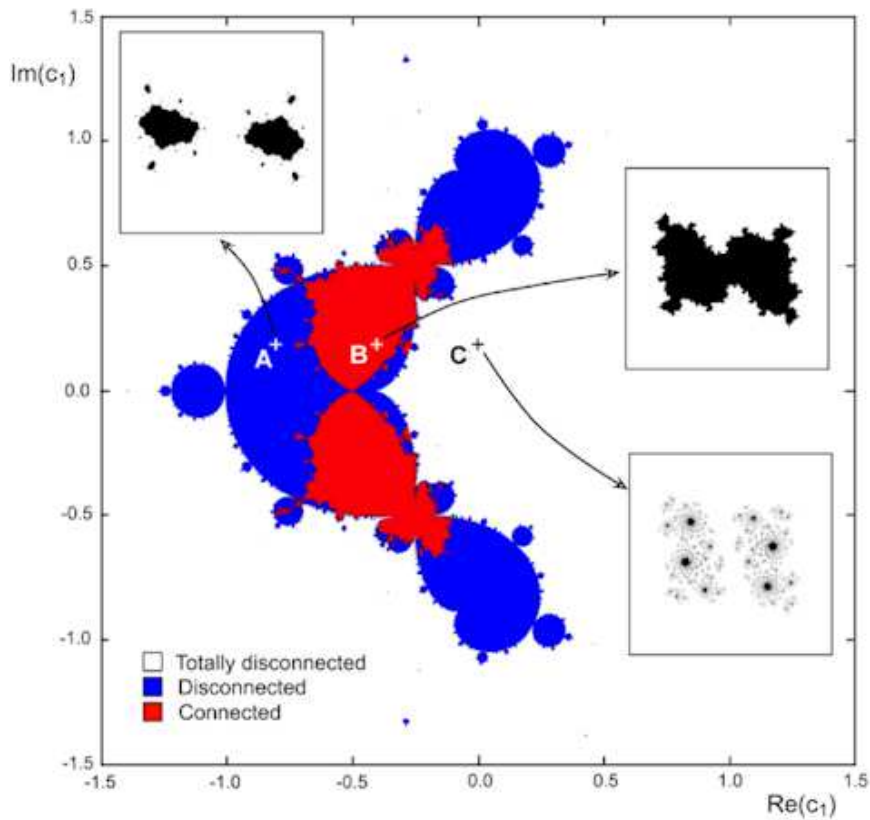}
\caption{Planar section obtained for $Re(c_2) = -0.4$, $Im(c_2) = 0$. The points $A(-0.8,0.2)$, $B(-0.4,0.2)$, and $C(0,0.2)$ respectively, belong to the disconnected, connected and totally disconnected respectively alternated Julia sets}
\label{fig4}       
\end{center}
\end{figure*}

\begin{figure}
\begin{center}
  \includegraphics[clip,width=0.55\textwidth] {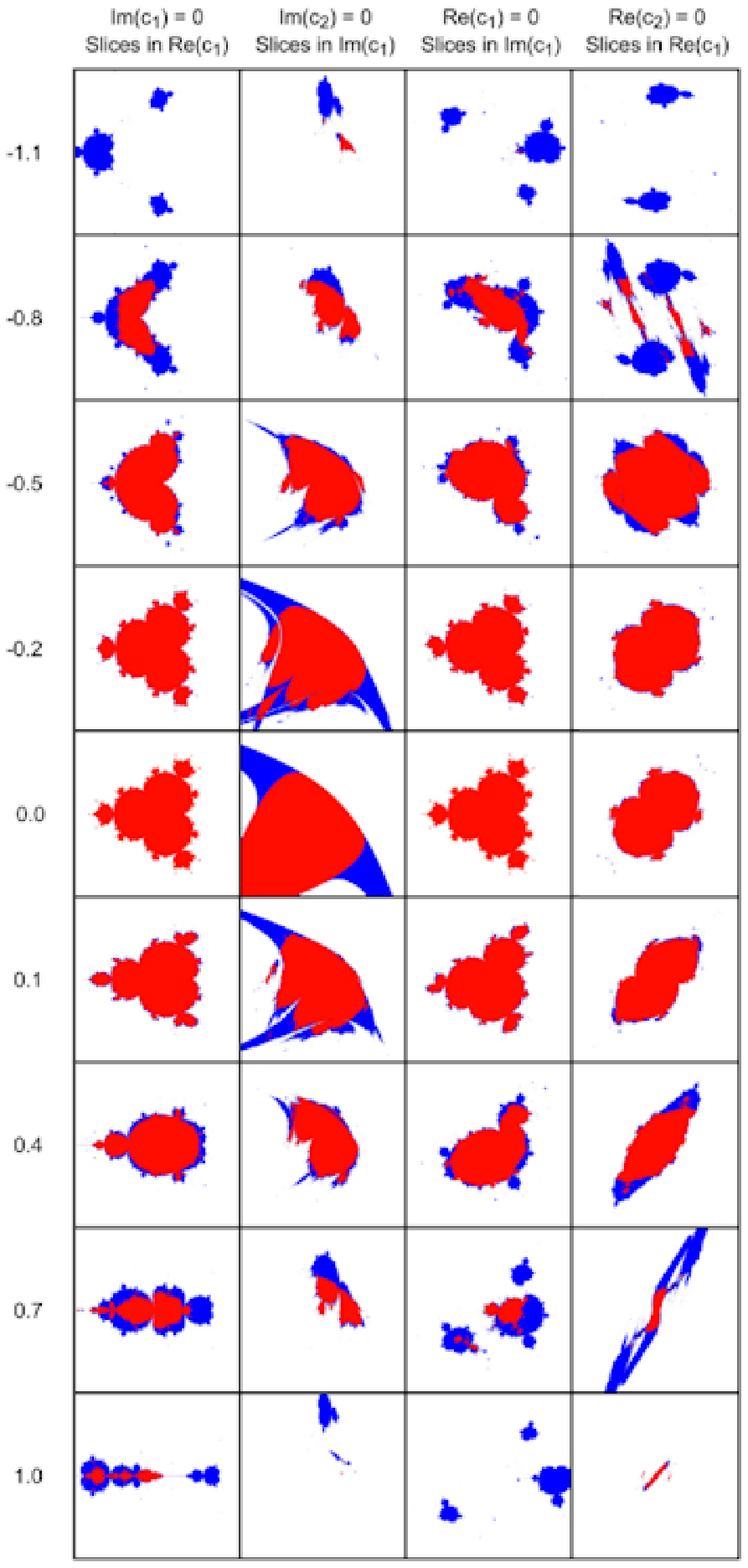}
\caption{Consecutive two-dimensional slices, through the 3D bodies $Im(c_1)=0$, $Im(c_2)=0$, $Re(c_1)=0$ and $Re(c_2)=0$, illustrating the transformation of the connectivity sets}
\label{fig5}       
\end{center}
\end{figure}

In 1992, Branner and Hubbard stated a conjecture for the general polynomial case of any degree \cite{bran1} that was subsequently proved by Branner and Hubbard \cite{bran2}. It completes the Fatou-Julia theorem as follows: The Julia set is totally disconnected if and only if each component of the filled Julia set containing critical points is aperiodic \cite{blanc,yy}.

	The disconnectedness means that the filled Julia sets consist of infinitely many pieces. Some of these pieces could be points, while others could be connected sets that are not points. In \cite{dan1} we proved that the filled Julia set of the alternated iteration of two quadratic family maps can be connected, totally disconnected and also disconnected (but not totally disconnected). The proof was made by employing analytical techniques.

	In this paper we draw and explore, with the aid of leading edge computer graphics algorithms, the fractal sets of the parameters which represent the three states: \emph{connectedness}, \emph{disconnectedness} and \emph{totally disconnectedness}, underlining the results previously obtained but also unveiling some new ones. Therefore, the obtained graphical results could represent possible subjects for future analytical works.

The loci or \emph{connectivity classes} (since each Julia set falls into exactly one of the three) corresponding to each of the three connectivity properties of alternated Julia sets, will be denoted by $CL$ (Connectedness Locus), $DL$ (Disconnectedness Locus) and $TDL$ (Totaly Disconnectedness Locus).

In order to graphically study these sets which, as we shall see are four dimensional, we employ a range of techniques for visualizing high dimensional objects more commonly reserved for studying 3D volumes arising from MRI or CT scans, techniques known as volume visualization.
	The structure of the paper is as follows: In Section 2 the alternated Julia sets and their connectivity properties are presented, while in Section 3 the connectedness is graphically analyzed.

\section{Connectivity sets of the alternated Julia sets}
\label{sec:1}

Let us consider the alternated Julia set generated by iterating the following system
\begin{eqnarray}\label{unu}
z_{n+1}=\left\{
\begin{array}{c}
z_{n}^{2}+c_{1},~~\text{if}~~n~~\text{is even,} \\
z_{n}^{2}+c_{2},~~\text{if}~~n~~\text{is odd.}%
\end{array}%
\right. ~~~\left( P_{c_{1}c_{2}}\right)
\end{eqnarray}
The alternated Julia set ${P_{{c_2}{c_1}}}$ can be defined in a similar way, i.e. $z_{n+1}$ takes the value $z_n^2 + {c_2}$ if $n$ is even and $z_n^2 + {c_1}$ if $n$ is odd. Because of the obvious symmetry of ${P_{{c_1}{c_2}}}$ and ${P_{{c_2}{c_1}}}$ we shall consider only the ${P_{{c_1}{c_2}}}$
case, the corresponding Julia set and filled Julia set being denoted ${J_{{P_{{c_1}{c_2}}}}}$ and ${K_{{P_{{c_1}{c_2}}}}}$ respectively. Moreover, because of Remark \ref{rem1}, we will study only $K$ sets.

Next, for the sake of simplicity hereafter, until necessary, we will drop the index ${c_1}{c_2}$.

The orbit generated by the iterates of $\{ P_{}^{ \circ n}\}$ is
\begin{eqnarray}
 z_0&, \\
 {z_1}& =& z_0^2 + {c_1},\\
 {z_2}& =& {(z_0^2 + {c_1})^2} + {c_2},\\
 {z_3}& =&{({(z_0^2 + {c_1})^2} + {c_2})^2} + {c_1},\\
 {z_4} & =& {({({(z_0^2 + {c_1})^2} + {c_2})^2} + {c_1})^2} + {c_2},\\
\vdots
\end{eqnarray}

Because the study of the dynamics of $P$ is a difficult task, the auxiliary quartic polynomial $Q$ shall be introduced
\begin{equation}\label{doi}
Q(w) = {({w^2} + {c_1})^2} + {c_2}.
\end{equation}

The (un)boundedness of $Q^{ \circ n}$ implies the (un)boundedness of $\{ {P^{ \circ n}}\} $ respectively \cite{dan1}. Thus, using the connectivity properties of $Q$ indicated by the general case of the Fatou-Julia theorem, it is proved that this theorem also
applies to $P$. To test numerically which values of the parameter $c$ give rise to connected, disconnected or totally disconnected filled Julia sets, without drawing the whole Julia set for each value of $c$, we have to test the boundedness of the orbits of the critical points.

According to the Fatou-Julia theorem, there are three possibilities for the orbits of the critical points of $Q$:

\begin{itemize}
\item [(i)]	The orbits of the critical points $0$ and $\pm \sqrt {-c_1}$ are bounded, the alternated Julia set $K$ being \emph{connected}.
\item [(ii)]The orbits of $0$ and $\pm \sqrt {-c_1}$ are unbounded, the alternated Julia set $K$ being\emph{ totally disconnected}.
\item [(iii)]Either the orbit of $0$ is bounded and the orbits of $\pm \sqrt {-c_1}$ are unbounded, or the orbit of $0$ is unbounded and the orbits of $\pm \sqrt {-c_1}$ are bounded. In this case, $K$ is \emph{disconnected} if the bounded orbits are periodic and \emph{totally disconnected} if the critical orbits are aperiodic.
\end{itemize}

\section{Visualization of the connectivity sets}

In this section we describe the way in which the connectivity domains are plotted.

Let us denote with $B$ the set of all complex numbers $c_1, c_2\in \mathbb{C}$ which determine the connectivity sets of alternated Julia sets. $B$ can be considered as a set of points of four coordinates $(Re(c_1),Im(c_1),Re(c_2),Im(c_2))$ in a 4D hyperspace verifying under the iteration (\ref{unu}) one of the relations (i), (ii) or (iii). In order visualize $B$, it can be imagined in the 4D - Quaternion space \cite{q}. Therefore, in the underlying base $\{1,i,j,k\}$, $B$ can be described as
\[B =  {{\mathop{\rm Re}\nolimits} ({c_1}) + {\mathop{\rm Im}\nolimits} ({c_1})i + {\mathop{\rm Re}\nolimits} ({c_2})j + {\mathop{\rm Im}\nolimits} ({c_2})k} ,\,\,\,{c_1},{c_2} \in  \mathbb{C},\]

\noindent where $i, j$, and $k$ are imaginary numbers satisfying the fundamental formula for quaternion multiplicative identities ${i^2} = {j^2} = {k^2} = ijk =- 1$, and forming the quaternion algebra.

Therefore, $B$ will be composed by three 4D combined sets corresponding to the three connectivity domains: $CL$, $DL$ and $TDL$.

By using volume rendering it is possible to study two-dimensional and three-dimensional sections through $B$. To create the 3D volume we fix one of the four coordinates: $\mathop{\rm Re}\nolimits ({c_1})$, $\mathop{\rm Im}\nolimits ({c_1})$, $\mathop{\rm Re}\nolimits ({c_2})$, or $\mathop{\rm Im}\nolimits ({c_2})$, this being equivalent to slicing the four-dimensional body $B$ with a hyperplane, for example $Re(c_1)=const$. Similarly, one obtains different three-dimensional sections by slicing with other hyperplanes, noting that these hyperplanes could be at arbitrary angles to the coordinate axes. If we fix two coordinates, for example ${\mathop{\rm Re}\nolimits} ({c_1})= {k_1}$ and ${\mathop{\rm Im}\nolimits} ({c_1}) = {k_2}$, with $k_1$ and $k_2$ some real constants (there are $6$ possible choices for two-dimensional slicesª $(Re(c_1),Im(c_1))$, $(Re(c_1),Im(c_2))$ and so on) one obtain two-dimensional images.

In this paper, in order to obtain the 3D bodies, taking into account the symmetries, we consider the following hyperplane slices: ${\mathop{\rm Re}\nolimits} ({c_1}) = 0$ (Figure \ref{fig1}a), ${\mathop{\rm Im}\nolimits} ({c_1}) = 0$, ${\mathop{\rm Re}\nolimits} ({c_2}) = 0$, and ${\mathop{\rm Im}\nolimits} ({c_1}) = 0$.

 Volume rendering requires a discrete, finite and generally regular sampling within a rectangularly bounded region. Therefore, in order to obtain an acceptable resolution of the main structures, the three-dimensional bodies were evaluated in a cubic region defined by the coordinates $(-2,-2,-2,)$ to $(2,2,2)$.

The volumes are sampled such that three different values are assigned to the three regions $CL$, $DL$ and $TDL$. Volume visualization involves forming a mapping between the values at each position in the regular sampling and colour/transparency. In the figures throughout this paper, the connected regions are assigned red, the disconnectedness regions blue, and the totally disconnectedness set (being a complement of the connectivity and disconnectivity sets in the considered 4D space) is made entirely transparent.

Volume rendering (see e.g. \cite{ren}) is the process by which it is determined how this volume, consisting of red and blue solid regions and other transparent regions, would appear if it were a real/physical object.

In Figure \ref{fig1}b an alternative visualization is achieved by projecting the connectivity set along the three orthogonal axes ${\rm{Im}}({c_1})$, ${\mathop{\rm Im}\nolimits} ({c_2})$ and ${\mathop{\rm Re}\nolimits} ({c_2})$ respectively.

In Figure \ref{fig2}, the disconnectedness regions corresponding to the cases ${\mathop{\rm Im}\nolimits} ({c_1}) = 0$, ${\mathop{\rm Re}\nolimits} ({c_2}) = 0$, and ${\mathop{\rm Im}\nolimits} ({c_2}) = 0$, are additionally made slightly transparent in order to reveal the internal connectedness regions. This internal structure of $B$ can additionally be revealed by taking planar slices through the volume, moving the cutting plane along one of the axis (see Figure \ref{fig3} where, again, the case ${\mathop{\rm Re}\nolimits} ({c_1}) = 0$  has been illustrated, the successive planes being orthogonal to the axis $Im(c_1)$).

If one considers the 2D slice obtained for $Re(c_2) = -0.4,~Im(c_2) = 0$ (Figure \ref{fig4}), to the points $A(-0.8,0.2)$, $B(-0.4,0.2)$, and $C(0,0.2)$, chosen from $DL$, $CL$ and $TDL$ respectively, generate disconnected, connected and totally disconnected alternated Julia sets.

Consecutive ``empirically'' chosen two-dimensional slices for other cases are presented in Figure \ref{fig5}, where each column represents slices, through the four 3D bodies, corresponding to $Im(c_{i})=0$ and $Re(c_{i})=0$, $i=1,2$ respectively.

\begin{remark}\label{rem2}
As can be seen in Fig. 1b, there are some separated ``parts'' and ``spots'' outside of the 3D bodies, but this ``disconnectedness'' phenomenon is probably due to the limited resolution at which the volume is sampled. Therefore, by using a higher resolution volume sampling (e.g. distance estimators for points near $\mathcal{M}$ set \cite{peit}) these apparently disconnected parts would be shown to be linked to the set.

\end{remark}

\section{Graphical properties of the connectivity sets of alternated Julia sets}

Analyzing the graphical results obtained in the previous section, a few properties can be emphasized such as:

\begin{itemize}
\item
The disconnectedness and totally disconnectedness sub-bodies are unbounded while the connectedness sub-bodies are bounded (see e.g. Fig. \ref{fig1}b and \ref{fig2});

\item
There exist typical symmetries like rotational symmetry, symmetry about origin or about one of the axis (see Fig. \ref{fig5});

\item There are planar sections, where $CL$ or $DL$ may have the shape of the $\mathcal{M}$ set, or Mandelbrot-like sets (``Multibrot'' sets) which correspond to the iterator $z^4+c$ (see e.g. the first column in Fig. \ref{fig5} for $Re(c_1)=0$).

\end{itemize}

\begin{figure}[ht]
\begin{center}
  \includegraphics[clip,width=0.4\textwidth]{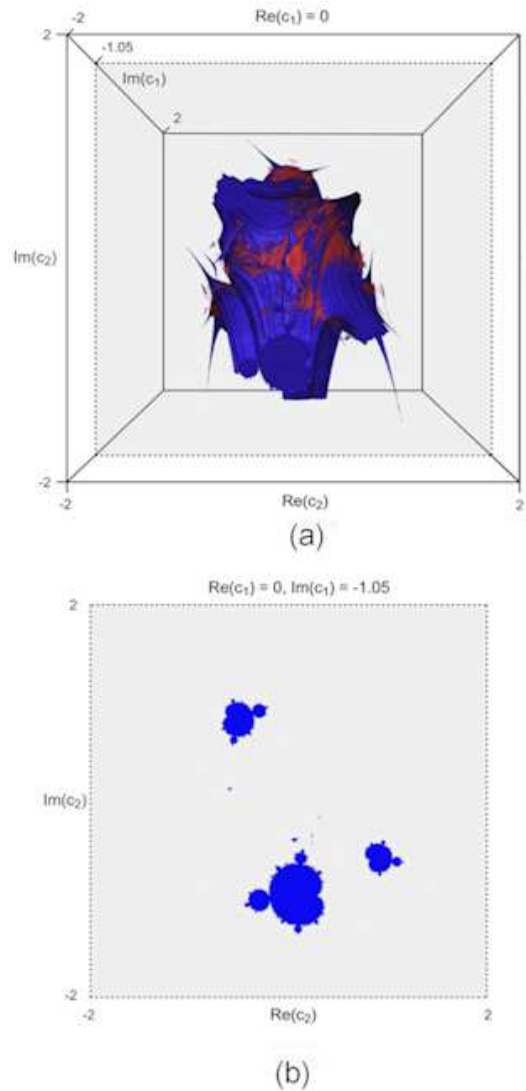}
\caption{$\mathcal{M}$ set found as cross sections of the disconnectivity sets. (a) Three-dimensional view. (b) Slice $Re(c_1)=0,~Im(c_1)=-1.05$ (Color figure online)}
\label{fig6}       
\end{center}
\end{figure}

\begin{figure}[h]
\begin{center}
  \includegraphics[clip,width=0.45\textwidth]{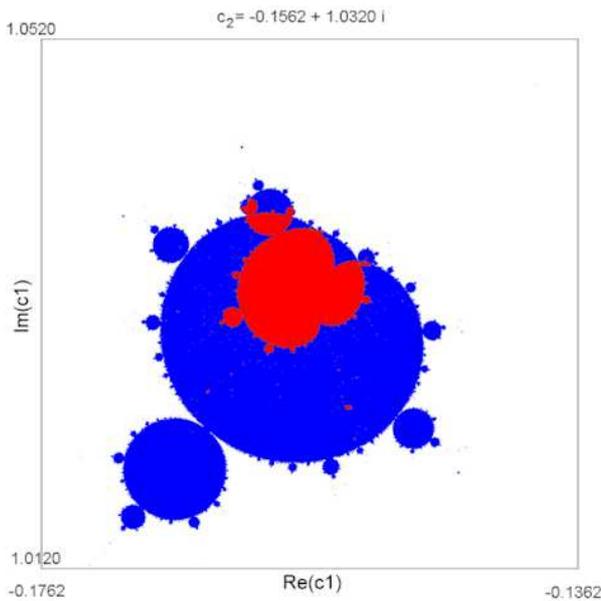}
\caption{$\mathcal{M}$-like sets in the slice $c_2=-0.1562+1.0320i$ (Color figure online)}
\label{fig7}       
\end{center}
\end{figure}

However, the most important property graphically unveiled in this paper, is a variant of one of the $\mathcal{M}$'s property, which is used as one of its definitions: $\mathcal{M}$ is the set of all values of the parameter $c$ for which each corresponding standard Julia set is connected. Beside this property, by locking to the ``ends'' of the three dimensional slices through $B$ in Fig. \ref{fig2}, the first and last rows in Fig. \ref{fig5}, or at the Fig. \ref{fig6},\ref{fig7}, we can deduce the following new property

\begin{property}
There are two-dimensional slices through $B$ where the set of all parameter values for which each corresponding alternated Julia set is
disconnected is the $\mathcal{M}$ set or, at least, a set resembling a Mandelbrot set.
\end{property}

Another example is presented in Fig. \ref{fig5}, where in the first image in the first column, obtained by sectioning the body corresponding to $Im(c_1)=0$, with  $c_1=-1.1+i0$, indicates that sets for which the alternated Julia sets are disconnected, are $\mathcal{M}$ sets. In the last image of the first column (corresponding to the same body), one can see that $\mathcal{M}$ corresponds either to connected alternated Julia sets (red color) or to disconnected Julia sets (blue) (see also Fig. \ref{fig7}).

\section{Conclusion}
In this paper we determined numerically and drew via computer graphics the connectivity loci of the alternated Julia sets introduced in \cite{dan1} and it appears that, the alternated Julia sets agree with a general form of Fatou-Julia Theorem. Thus, the alternated Julia sets can be connected, disconnected or totally disconnected.

Moreover, we verified visually some other new properties, the most important being that the $\mathcal{M}$ set represents the set of all points in the complex plane for which the alternated Julia sets are disconnected (but not totally disconnected).

Therefore, even though this work does not substitute for the underlying proof, it could represent a guide for future research.

\vspace{5mm}
\textbf{Acknowledgements}

We thank Robert L. Devaney and Bodil Branner for their useful comments which helped us
to understand better the connectivity aspects presented in \cite{dan1}.

We thank the reviewers for their valuable comments which helped to considerably improve the quality of the manuscript.

The work was supported by iVEC through the use of advanced computing resources located at the University of Western
Australia. The volume rendering is performed using Drishti, a ``Volume Exploration and Presentation Tool'' developed by Ajay Limaye at the Australia National University.

%

\end{document}